\documentclass[11pt,a4paper,reqno]{article}

\usepackage{graphicx}
\usepackage{amssymb}
\usepackage{verbatim}
\usepackage[cp1250]{inputenc}
\usepackage[T1]{fontenc}
\usepackage{graphics}
\usepackage{color}
\usepackage[color,matrix,arrow]{xy}
\usepackage[all,color]{xypic}
\usepackage{hyperref}

\newtheorem {definition}{Definition}
\newtheorem{theorem}{Theorem}
\newtheorem {bemerkung}{Remark}
\newtheorem{beispiel}{Example}
\newtheorem{frage}{Question}
\newtheorem{vermutung}{Conjecture}

\newenvironment{remark} {\begin{bemerkung} \normalfont }{\end{bemerkung}}
\newenvironment{example} {\begin{beispiel} \normalfont }{\end{beispiel}}

\def\ZZ{\mathbb{Z}}
\def\MM{\mathbb{M}}
\def\NN{\mathbb{N}}
\def\DD{\mathbb{D}}

\date{}

\begin{document}
\title{On tame strongly simply connected algebras}

\author{Stanisław Kasjan and Andrzej Skowroński}

\maketitle

\thanks{\em Dedicated to Jos\'e Antonio de la Pe\~na on the occasion of his 60th birthday.}

\begin{abstract}{In this survey we present the criterion for tameness of strongly simply connected algebras due to Br\"ustle, de la Pe\~na and Skowroński. We recall relevant concepts of representation theory and discuss some applications and connections to other problems.}\end{abstract}

\section{Introduction}

The criterion of tameness of strongly simply connected algebras - Theorem \ref{main} below, due to Thomas Br\"ustle, Jos\'e Antonio de la Pe\~na and Andrzej Skowroński \cite{BrPeSk2011} can be considered as a culmination of a research program aimed on
better understanding of the concept of representation types of algebras. The purpose of this note is to present the result and the on a background of classical concepts and results of representation theory, indicating also some less obvious connections to other mathematical ideas. For more information the reader is referred to \cite{PeSk2011}.

{\begin{theorem}\cite[Main Theorem]{BrPeSk2011} \label{main}{ Let $A$ be a strongly simply connected algebra. The following statements are equivalent:
\begin{enumerate}
\item $A$ is a tame algebra.
\item The Tits form $q_A$ of $A$ is weakly nonnegative.
\end{enumerate}}\end{theorem}

Throughout this article $A$ is an associative algebra of finite dimension over an  algebraically closed field $K$.
Usually we assume that $A$ is basic, equivalently, by classical Gabriel's theorem, $A$ is isomorphic to a bound quiver algebra, that is

 $$A\cong KQ/I,$$
for some finite quiver $Q$ and an admissible ideal $I$ in the path algebra $KQ$. We refer to \cite{AsSiSk2006} for more information about the basic concepts of representation theory of bound quivers and algebras.

 The general question of the representation theory of finite dimensional algebras can be formulated as follows:
 classify (indecomposable) right $A$-modules of finite $K$-dimension up to isomorphism. Of course, we are far from satisfactory solution of the problem raised this way.

Finite dimensional modules can be represented  by sequences of matrices and this can be done  in several ways. One way is to consider representations of the bound quiver $(Q,I)$ associated to the algebra $A$, another - based on the idea of analyzing  module through its projective presentation - leads to the concept of Roiter's matrix problems \cite{Ro1980} or  BOCS introduced by Drozd \cite{Dr1972}.

Finding isomorphism classes of $A$-modules is equivalent to  finding "canonical forms" of matrices up to certain admissible operations. Let us remark that matrix problem algorithms allow to
 obtain a classification of indecomposable $A$-modules of  $K$-dimension bounded by a number fixed \'a priori. That does not mean, of course,  that a complete classification is always possible.

Let us recall some important facts concerning the easiest situation - the case of
representation-finite algebras. By the
definition, an algebra
 $A$ is {\em representation-finite} (or {\em of finite representation type}) if it allows only finitely many isomorphism classes of finite dimensional indecomposable $A$-modules. By famous
{\em 1st Brauer-Thrall Conjecture}, which is now a theorem by   Roiter \cite{Ro1968} and Auslander \cite{Au1978}, $A$ is representation-finite if and only if there is a common bound $m_A$ for the dimensions of indecomposable $A$-modules. Moreover, if $A$ is not representation-finite (remember that $K$ is algebraically closed), then there exist infinitely many\footnote{By \cite{Sm1980} it is enough to assume that there exists one such $m$.}  numbers $m\in\NN$ such that for each of those $m$ there is infinitely many pairwise non-isomorphic $A$-modules of $K$-dimension $m$. This is the content of the {\em 2nd Brauer-Thrall Conjecture} confirmed by several authors, see  \cite{Ba1985},\cite{Bo1985}, \cite{BrTo1986},  \cite{Fi1985}, (also \cite{NaRo1973}).
The reader is referred to
\cite{Ri1979} for more on Brauer-Thrall conjectures. Validity of the 2nd Brauer-Thrall conjecture allows to prove that given a natural number $d$, the class of $d$-dimensional representation-finite algebras over algebraically closed fields is finitely axiomatizable in a suitable first order language, see
 \cite{JeLe1981}, \cite{JeLe1982}, \cite[Theorem 12.54]{JeLe1989}.

This corresponds to the fact, that given na natural number $d$ there is only finitely many isomorphism classes of representation-finite algebras having the $K$-dimension equal $d$. This is, in turn, a consequence of another deep and difficult result obtained in the 80's of the 20th century: the proof of the existence of a multiplicative basis of every representation-finite algebra over an algebraically closed field \cite{BaGaRoSa1985}.

 Note that Jensen and Lenzing's arguments for the finite axiomatizability are not constructive.  However, this can be made more explicite by applying a "numerical version" of 2nd Brauer-Thrall Conjecture, see \cite[2.4]{Ge1995}. Finally, it is possible to verify in a finite number of steps
 if  given an algebra $A$  is representation-finite or not. In practice, for this purpose we use Bongartz's criterion for representation-finiteness, which shall be discussed later, Auslander-Reiten theory, Galois covering techniques etc.  These methods allow, in the representation-finite case, to give a list of representatives of the isomorphism classes of the indecomposable $A$-modules.

\medskip

If $A$ is not representation-finite, then there are infinitely many isomorphism classes of indecomposable modules of dimension $m$ for infinitely many numbers $m$. It may happen that for every $m$, all but finitely many representatives of those classes can be arranged in finitely many "rational 1-parameter families". In this case we say that $A$ is {\em of tame representation type} or simply {\em tame}\footnote{ According to our convention a representation-finite algebra is tame.}. We omit precise formulations of the well-known concepts appearing here, the reader is referred to a reach literature on the subject, see for example \cite{Ri1984}, \cite{Si1992}, \cite{SiSk2007}.

 Otherwise the problem of classification of two non-commuting square matrices up to simultaneous conjugations can be "embedded" into the problem of classifying the indecomposable $A$-modules, in which case we say that $A$ is {\em of wild representation type}  or {\em wild}.

 This is essentially the content of the famous {\em Tame-Wild Dichotomy}: \medskip

\begin{theorem} \cite{Dr1979}, \cite{CB1988} { $A$ is either tame or  wild.}\end{theorem}\medskip

We can not hope obtain a complete classification of isomorphism classes of indecomposables in the wild case. Also in the tame case this task can be unreachable. The level of complication is reflected somehow by the stratification of the class of tame algebras into subclasses: domestic, polynomial growth. Let us recall the notions introduced  in  \cite{Sk1987}, see also  \cite{Sk1990}. See also \cite[Section 14.4]{Si1992}, \cite[Chapter XIX]{SiSk2007}.

Let $A$ be a tame algebra. Given a natural number $m$ let $\mu_A(m)$ be the lest number of 1-parameter families needed to parameterize (all but finite number) indecomposables of dimension $m$.
\begin{itemize}
\item If $\mu_A(m)$ is bounded, then $A$ is {\em finite growth}. By the results of \cite{CB1991} an algebra is of finite growth if and only if it is {\em domestic} in the sense of \cite{Ri1980}.

\item If $\mu_A(m)\le m^N$ for some $N$ and arbitrary $m$, then $A$ is said to be of {\em polynomial growth}.
\end{itemize}

We have remarked that given an algebra $A$ it is possible to verify wether $A$ is representation-finite or not and it is related with finite axiomatizability of representation-finiteness. This fact is also connected with Gabriel's result that "finite representation type is open" \cite{Ga1974}, more precisely:  representation-finite algebras of fixed dimension $d$ induce Zariski-open subset in the variety of $d$-dimensional associative algebras. Indeed, an analysis of the proof of Geiss's result on degenerations of wild algebras \cite{Ge1995} shows that the openness of the varieties of representation-finite algebras can be derived from the existence of finite number of axioms of representation-finiteness.

Let us collect some known analogies for tameness. First of all,
the class of tame algebras (of fixed dimension, over algebraically closed fields) is axiomatizable \cite{Ka2002}. This is related with the observation, that if $A$ is wild, then we can verify it in finite number of steps\footnote{By applying, for example, matrix problems techniques: we classify canonical forms of indecomposable modules of dimension $m$,  starting from $m=1$ and increasing $m$ until we meet 2-parameter family of pairwise non-isomorphic indecomposable $m$-dimensional modules.}.
However we do not know if tame is finitely axiomatizable, that means that it is possible, that there is an algebra $A$ such that we will never be sure if $A$ is tame or not.
If we knew that tame is finite axiomatizable, we would prove that "tame is open"  \cite{Ka2002}, by the same kind of arguments that we mentioned above with respect to openness of finite representation type.

Theorem \ref{main} provides us with a criterion for tameness valid for strongly simply connected algebras. The condition equivalent to tameness is expressed  in term of Tits quadratic form.

\section{Tits quadratic form}

Let $A$ be a basic finite dimensional $K$-algebra of the form $KQ/I$, where
 $Q=(Q_0,Q_1,s,t)$ is a finite acyclic quiver\footnote{As usual, $Q_0$ is the set of vertices, $Q_1$ is the set of arrow and $s(\alpha)$ and $t(\alpha)$ denote the source  and the target of the arrow $\alpha$.} and
$I$ is an admissible\footnote{That means $(KQ_1)^m\subseteq I\subseteq (KQ_1)^2$ for some $m$.} ideal of the path algebra $KQ$. Recall that by a {\em relation} in the path algebra $KQ$ we mean a linear combination of paths having the same sink and the same source.  Let $R$ be a minimal set of relations generating the ideal $I$. $R$ decomposes into the disjoint union
$\bigcup_{i,j\in Q_0}R_{i,j}$, where $R_{i,j}$ denote the set of the elements starting at the vertex $i$ and ending at $j$. Although $R$ is not uniquely determined by $I$, the numbers
$r_{i,j}:=|R_{i,j}|$ do not depend on the set $R$ \cite{Bo1983}.

The {\em Tits quadratic form} of the algebra $A$ is the integer quadratic form
$$q_A:\ZZ^{Q_0}\rightarrow \ZZ$$
defined by the formula
\begin{equation}\label{tits} q_A(z)=\sum_{i\in Q_0}z^2_i-\sum_{i\rightarrow j\in Q_1}z_iz_j+\sum_{i,j\in Q_0}r_{i,j}z_iz_j
\end{equation}
for $z=(z_i)_{i\in Q_0}\in\ZZ^{Q_0}$.

Let us also recall another quadratic form associated with an algebra: the
{\em Euler quadratic form}. For this purpose denote by
$S_i$ the simple $A$-module corresponding to the vertex $i\in Q_0$ and assume that $A$ has finite global dimension. The quadratic form

$$
\chi_A:\ZZ^{Q_0}\rightarrow \ZZ,
$$
defined by the formula
\begin{equation}\label{euler}
\chi_A(z)=\sum_{s=0}^{\infty}\sum_{i,j\in Q_0}(-1)^s\dim{\rm Ext}_A^s(S_i,S_j)z_iz_j.
\end{equation}
is called the  {\em Euler quadratic form} of the algebra $A$.

The Euler form is often called a homological form, by obvious reasons, whereas the Tits form - the geometric quadratic form of an algebra. This is because the summands of the right hand side of (\ref{tits}) have interpretations as (bounds of) dimensions of some geometric objects associated with $A$ and $z$. Thanks to this interpretation one can prove that:

\begin{enumerate}
\item If $A$ is representation-finite, then $q_A$ is weakly positive, that is, $q_A(z)>0$ for all non-zero vectors $z$ with all coordinates nonnegative \cite{Bo1983}.

\item If $A$ is tame, then $q_A$ is weakly nonnegative, that is, $q_A(z)\ge 0$ for all vectors $z$ with all coordinates nonnegative \cite{Pe1991}.
\end{enumerate}

The assertion 1. is the easier part of the Bongartz's criterion for representation-finiteness, we shall discuss the converse later. The converse of 2. (valid under suitable conditions on $A$) is the content of Theorem \ref{main}. Let us recall the main steps of the proof of 2, following de la Pe\~na.

It is well-known that the $A$-modules of dimension vector $z=(z_i)_{i\in Q_0}$ can be represented by the tuples $(M_{\alpha})_{\alpha\in Q_1}\in\prod\MM_{z_{t(\alpha)}\times z_{s(\alpha)}}(K)$ of matrices satisfying the relations defining the ideal $I$. These tuples form the {\em variety of modules} of dimension vector $z$ denoted by ${\rm mod}_A(z)$ and, by Krull Theorem, it is easy to observe that
\begin{equation}\label{weaknon1}
\dim{\rm mod}_A(z)\ge \sum\limits_{i\rightarrow j\in Q_1}z_iz_j-\sum\limits_{i,j\in Q_0}r_{i,j}z_iz_j.
\end{equation}

The isomorphism classes of modules of dimension vector $z=(z_i)_{i\in Q_0}$ can be identified with the orbits of the natural algebraic group action
$$
Gl(z)\times {\rm mod}_A(z)\rightarrow {\rm mod}_A(z)
$$
of the group $Gl(z)=\prod_{i\in Q_0}Gl_{z_i}(K)$. Clearly
\begin{equation}\label{weaknon2}
\dim Gl(z)=\sum_{i\in Q_0}z_i^2.
\end{equation}
Now, by (\ref{tits}), (\ref{weaknon1}) and (\ref{weaknon2}),
\begin{equation}\label{weaknon3}
q_A(z)\ge \dim Gl(z)-\dim {\rm mod}_A(z).
\end{equation}

Now it is time to apply de la Pe\~na's argument: if $A$ is tame,  representatives of the indecomposables of fixed dimension can be arranged into finite number of "one-parameter families". Therefore, there is a "$|z|$-parameter family" of  representatives of the modules of dimension vector $z$, where $|z|=\sum\limits_{i\in Q_0}z_i$. It follows that if
 $A$ is tame, then
 \begin{equation}\label{weaknon4}
 \dim Gl(z)+ \sum\limits_{i\in Q_0}z_i\ge \dim{\rm mod}_A(z).
 \end{equation}

 Combining (\ref{weaknon3}) and (\ref{weaknon4}) we obtain
 $$q_A(z)\ge \dim Gl(z)-\dim {\rm mod}_A(z)\ge - \sum\limits_{i\in Q_0}z_i.$$
As the left hand side of the above inequality depends quadratically on $z$, whereas the right hand side - linearly, we conclude that
$q_A(z)\ge 0$ for every $z$ with nonnegative coordinates and the proof is complete.

As it is mentioned above, the converses of 1. and 2. are not true in general. Let us recall the well known Bongartz's example showing this.

\begin{example}\label{Bongartzex} Let $Q$ be the quiver
{\small
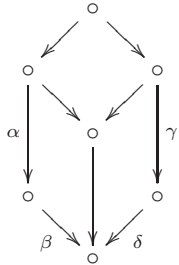
$$\xymatrix@=0.3pc{&&\ar[ddll]\circ\ar[ddrr]&&\\
&&&&\\
\circ{\ar@[black][dddd]_{\alpha}}\ar[ddrr]&&&&\circ{\ar@[black][dddd]^{\gamma}}\ar[ddll]\\
&&&&\\
&&\circ\ar[dddd]&&\\
&&&&\\
\circ{\ar@[black][ddrr]_{\beta}}&&&&\circ{\ar@[black][ddll]^{\delta}}\\
&&&&\\
&&\circ&&}
$$}
\end{example}
Assume that $I_1$ is the ideal  generated by all commutativity relations, and $I_2$ the ideal is generated by the commutativity of the upper square and two zero-relations $\alpha\beta$ and $\gamma\delta$. Then $A_1=KQ/I_1$ is representation-finite, $A_2=KQ/I_2$ is wild, and the Tits quadratic forms of $A_1$ and $A_2$ coincide and are weakly positive.

A lot of effort has been put into proving the converses of 1. and 2. under various assumptions. The most satisfactory results in this direction are the Bongartz's criterion for representation-finiteness and Theorem \ref{main} for tameness. The results can be considered as a culmination of a research programm scheduled by Sheila Brenner, who wrote in 1975 about her study on connections between the representation type and definiteness of certain quadratic forms \cite{Br1975}:
{ \em "This paper is written in the spirit of experimental science. It reports some observed regularities and suggests that there should be a theory to explain them."}

\section{Simply connected and strongly simply connected algebras}

Let $A=KQ/I$ be as before, assume that the quiver $Q$ is connected and let
$m(I)$ be a set of  minimal relations generating $I$.
A relation $\sum_{i\in I}\lambda_iu_i\in I$, where $\lambda_i\in K$ and $u_i$ are paths in $Q$,  is {\em minimal} if  $\sum_{i\in J}\lambda_iu_i\notin I$ for any proper subset $J\subseteq I$.
We look at
 $Q$ as a 1-dimensional complex and,
whenever two paths $u,v$ are involved into a minimal relation, we attach a 2-cell along the loop $uv^{-1}$:
$$
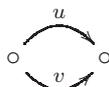
\xymatrix@=0.3pc{\circ\ar@/^1pc/[rrr]^u\ar@/_1pc/[rrr]^v&&&\circ}\medskip
$$
The fundamental group of the obtained complex is called the fundamental group of the bound quiver $(Q,I)$ and it is denoted by  $\pi_1(Q,I)$ \cite{MaPe1983}. Let us remark that $\pi_1(Q,I)$ depends not only on $KQ/I$,  it may depend on the particular choice of $I$.

\begin{definition} \cite{AsSk1988}.
The algebra $A$ is {\em simply connected} if $\pi_1(Q,I)$ is trivial for every presentation $A\cong KQ/I$.
\end{definition}

In the Example \ref{Bongartzex} we have $\pi_1(Q,I_1)$ trivial, whereas $\pi_1(Q,I_2)\cong\ZZ* \ZZ$ is a free non-commutative group with 2 free generators.

Combining the results of \cite{Bo1983}, \cite{Bo1984},   \cite{BaGaRoSa1985} and \cite{BrGa1984} we derive

\begin{theorem}\label{Bongartz}{ If $A$ is {\bf simply connected}, then $A$ is representation-finite if and only if  $q_A$ is weakly positive.}\end{theorem}

In fact, the theorem is valid for algebras having simply connected Galois covering.  Results on preservation of (locally) representation-finiteness by Galois coverings, \cite{BoGa1981}, are important ingredients in the proof of this general version.

The last sentence of the paper \cite{Bo1983}
is
{\em "There should be a similar result for tame representation type."}

Of course, it is expected that weak nonnegativity of the Tits form implies tameness under suitable assumptions. Actually, there are
 several results of that kind. Such implication has been proved for instance for:

\begin{itemize}
\item hereditary algebras \cite{Na1973},
\item tilted algebras \cite{Ke1989},
\item quasitilted algebras\cite{Sk1998}.
\end{itemize}

Observe that it is not enough to assume that $A$ is simply connected, as the following example shows:
\begin{example}\label{BPS}
Let $Q$ be the quiver
$$
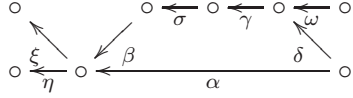
\xymatrix@=0.3pc{\circ&&&&\circ\ar[ddll]^{\beta}&&\circ\ar[ll]^{\sigma}&&\circ\ar[ll]^{\gamma}&&\circ\ar[ll]^{\omega}\\
&&&&&&&&&&\\
\circ&&\circ\ar[ll]^{\eta}\ar[uull]^{\xi}&&&&&&&&\circ\ar[llllllll]^{\alpha}\ar[uull]^{\delta}
}
$$
and let $I$ be the ideal generated by the relations:
$$\alpha\xi, \omega\gamma\sigma\beta, \alpha\eta-\delta\gamma\sigma\beta\eta$$

Then $A=KQ/I$ is
simply connected, the Tits quadratic form $q_A$ is weakly nonnegative, but $A$ is
  wild \cite[Example 1.7]{BrPeSk2011}.

\end{example}

Thus a strengthening of the concept of simply connectedness is needed.


\begin{definition}\label{ssc} \cite{Sk1992}
The $A$ is {\bf strongly simply connected} if  every convex subcategory $C$ of $A$ is simply connected.
\end{definition}

It is proved in \cite{Sk1992} that the condition in the definition above is equivalent to the vanishing of the first Hochschild cohomology groups $H^1(C)$ (with coefficients in the bimodule ${}_CC_C$) for every convex subcategory $C$ of $A$. By results of Bretscher - Gabriel \cite{BrGa1984}  and Martinez - de la Pe\~na \cite{MaPe1983}
 a triangular representation-finite algebra  is simply connected if and only if it is strongly simply connected.

\section{Strongly simply connected of polynomial growth}

The proof of Theorem \ref{main} depends essentially on some older result on polynomial growth algebras.

{\begin{theorem} \label{polygr}   { Let $A$ be a strongly simply connected algebra. The following statements are equivalent:
\begin{enumerate}
\item $A$ is  tame of polynomial growth.
\item The Tits form $q_A$ of $A$ is weakly nonnegative and $A$ does not contain a convex subcategory, which is $pg$-critical. .
\end{enumerate}
}\end{theorem}}

The list of $pg$-critical algebras is given in \cite{NoSk1996}.

Representation theoretical properties of tame polynomial growth algebras differ substantially from those of remaining tame algebras. In particular, we observe characteristic behavior of the Tits and Euler quadratic forms in case of polynomial growth strongly simply connected algebras. Let us recall two theorems:

{\begin{theorem} \cite{PeSk1996}  Let $A$ be a strongly simply connected algebra. The following conditions are equivalent:
\begin{enumerate}
{\item $A$ is tame of polynomial growth.}
{\item For every indecomposable module $X$ of dimension vector $z$
$$
\chi_A(z)=\dim G(z) - \dim {\rm mod}_A(z)\ge 0.
$$}
\end{enumerate}
{In this case every indecomposable module is a smooth point of the module scheme\footnote{Let us remark, that here we mean the module scheme, see \cite{PeSk1996} for the details, not the module variety ${\rm mod}_A(z)$}.}
\end{theorem}}\medskip

{\begin{theorem}\cite{PeSk1999} Let $A$ be a strongly simply connected algebra and $X$  an indecomposable $A$-module of dimension vector $z$. Then:
\begin{enumerate}
{\item $0\le q_A(z)\le 18|Q_0|$.}
{\item $0\le\chi_A(z)\le 2+\#\{\text{projective-injective indecomposable $A$-modules}\}$.}
{\item If $X$ is a faithful module, then $q_A(z)\le 2$ and $\chi_A(z)\le 2$.}
\end{enumerate}
\end{theorem}}\medskip

\begin{remark} There are also interesting connections between the
growth of tame algebra and certain phenomena on the level of infinite-dimensional modules.

 Let us recall that a module is {\em superdecomposable} if it has no indecomposable direct summands.
Here we assume that the field $K$ is countable.

\begin{enumerate}
{\item A wild algebra possesses a superdecomposable pure-injective module  \cite{GrPr2016}.}
{\item A non-domestic string algebra  possesses a superdecomposable pure-injective module \cite{Pu2004}.}
{\item An algebra admitting a strongly simply connected Galois covering of non-polynomial growth, in particular, a strongly simply connected algebra of non-polynomial growth, possesses a superdecomposable pure-injective module, under the additional assumption that  $char(K)\neq 2$ \cite{KaPa2014}, \cite{KaPa2016}.}
{\item Tubular algebras possess a superdecomposable pure-injective module \cite{HaPr2015}.}
{\item A strongly simply connected algebra is domestic if and only if does not possess a superdecomposable pure-injective module, under the additional assumption that  $char(K)\neq 2$ \cite{Pa2015}.}

\end{enumerate}
\end{remark}

\section{On the proof of the main theorem}

We know that the Tits quadratic form of a tame algebra is weakly nonnegative, so we concentrate on the converse implication for strongly simply connected algebras. Very briefly the proof can be sketched as follows.
By Theorem \ref{polygr} a strongly simply connected algebra $A$ is tame of polynomial growth if and only if the Tits quadratic form  $q_A$ is weakly nonnegative and $A$ does not contain a $pg$-critical convex subcategory.
 Thus we can restrict our attention  to $A=KQ/I$ with $q_A$ weakly nonnegative and  containing a $pg$-critical convex  subcategory. Without loss of generality we can assume that additionally $A$ has an indecomposable module whose support contains all sources and all sinks of $Q$. The structure of such algebras can be described quite precisely; they are so called
 $\DD$-algebras in the terminology of \cite{BrPeSk2011}. Moreover, to $A$ there is associated a {\em mild and smooth} algebra $A^*$, which is tame if and only if $A$ is tame. The reader is again referred to \cite{BrPeSk2011} for the details. So we need to prove that $A^*$ is tame. This final step is done by  Geiss's degeneration theorem \cite{Ge1995}, since
 $A^*$ degenerates to a special biserial algebra, which is tame  by \cite{WaWa1985},   \cite{DoSk1987},  \cite{BuRi1987}.

Let us show on examples ideas how to perform degenerations mentioned above. There are two basic, rather known,  tricks.
First is to degenerate a commutative square to a square with a zero relation:

$$\begin{array}{ccc}
\xymatrix@=0.3pc{&&\circ\ar[ddll]_{\alpha}\ar@{--}[dddd]\ar[ddrr]^{\beta}&&\\
&&&&\\
\circ\ar[ddrr]_{\gamma}&&&&\circ\ar[ddll]^{\delta}\\
&&&&\\
&&\circ&&}
&
\xymatrix@=0.3pc{&&\\&& \\ \ar@{~>}[rr]&&\\&&\\&&}
&
\xymatrix@=0.3pc{&&\circ\ar@[black][ddll]_{\alpha}\ar[ddrr]^{\beta}&&\\
&&&&\\
\circ\ar@[black][ddrr]_{\gamma}&&&&\circ\ar[ddll]^{\delta}\\
&&&&\\
&&\circ&&}\vspace{5ex}\\
\alpha\gamma=\beta\delta&&\alpha\gamma=0
\end{array}$$

The degenration is performed by the following "passage to the limit":
$\alpha\gamma-t\beta\delta \rightarrow \alpha\gamma$ as $t\rightarrow 0$.

Another trick is the following:

$$\begin{array}{ccc}
\xymatrix@=0.3pc{&&\circ\ar[ddll]_{\alpha}\ar@{--}[dddd]\ar[ddrr]^{\beta}&&\\
&&&&\\
\circ\ar[ddrr]_{\gamma}&&&&\circ\ar[ddll]^{\delta}\\
&&&&\\
&&\circ&&}
&\xymatrix@=0.3pc{&&\\&& \\ \ar@{~>}[rr]&&\\&&\\&&}
&
\xymatrix@=0.3pc{&&\circ\ar@[black][dd]_{{\alpha}}&&\\
&&&&\\
&&\circ\ar@[black]@(ul,dl)_{{\epsilon}}\ar@[black][dd]_{{\gamma}}&&\\
&&&&\\
&&\circ&&}\vspace{5ex}\\
\alpha\gamma=\beta\delta&&\alpha\gamma=0, \epsilon^2=0
\end{array}$$

This is based on the observation that the path algebra of the two point quiver without arrows ($\begin{array}{cc}\circ&\circ\end{array})$ is isomorphic to $K\times K \cong K[X]/(X(X-t))$ for $t\neq 0$ and $K[X]/(X(X-t)) \rightarrow K(\xymatrix@=0.3pc{\circ\ar@(ul,dl)_{{\epsilon}}})/\epsilon^2$ as $t\rightarrow 0$.

Combining this tricks we can degenerate the pg-critical algebra defined by the bound quiver

{
{\small
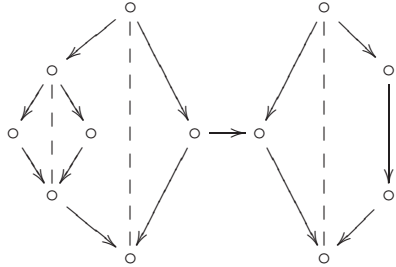
$$\xymatrix@=0.3pc{
&&&\circ\ar[ddll]\ar@{--}[dddddddd]\ar[ddddrr]&&&&&&\circ\ar[ddddll]\ar@{--}[dddddddd]\ar[ddrr]&&\\
&&&&&&&&&&&\\
&\circ\ar[ddl]\ar@{--}[dddd]\ar[ddr]&&&&&&&&&&\circ\ar[dddd]\\
&&&&&&&&&&&\\
\circ\ar[ddr]&&\circ\ar[ddl]&&&\circ\ar[rr]\ar[ddddll]&&\circ\ar[ddddrr]&&&&\\
&&&&&&&&&&&\\
&\circ\ar[ddrr]&&&&&&&&&&\circ\ar[ddll]\\
&&&&&&&&&&&\\
&&&\circ&&&&&&\circ&&}
$$}

to the special biserial algebra defined by the quiver

{\small
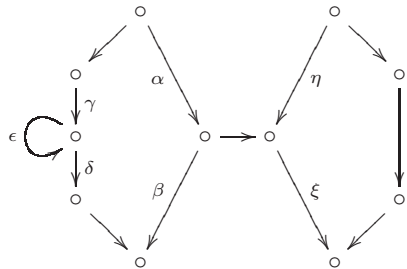
$$\xymatrix@=0.3pc{
&&&\circ\ar[ddll]\ar@[black][ddddrr]_{\alpha}&&&&&&\circ\ar@[black][ddddll]^{\eta}\ar[ddrr]&&\\
&&&&&&&&&&&\\
&\circ\ar@[black][dd]^{\gamma}&&&&&&&&&&\circ\ar[dddd]\\
&&&&&&&&&&&\\
&\circ\ar@[black][dd]^{\delta}\ar@[black]@(ul,dl)_{{\epsilon}}&&&&\circ\ar[rr]\ar@[black][ddddll]_{\beta}&&\circ\ar@[black][ddddrr]^{\xi}&&&&\\
&&&&&&&&&&&\\
&\circ\ar[ddrr]&&&&&&&&&&\circ\ar[ddll]\\
&&&&&&&&&&&\\
&&&\circ&&&&&&\circ&&}
$$}

with zero-relations $\alpha\beta=\gamma\delta=\eta\xi=\epsilon^2=0$.}

\section{Some consequences and comments}

The main theorem yields a possibility of checking in a finite number of steps wether given a strongly simply connected algebra is tame or not.  It is not difficult to observe that also it is possible to check algorithmically wether  $A=KQ/I$ is strongly  simply connected - here the Hochschild cohomology approach is useful. Thanks to this we obtain

{\begin{theorem} \cite{Ka2007}{  Given a number $d$, the class of $d$-dimensional tame strongly simply connected algebras over algebraically closed fields is finitely axiomatizable.}\end{theorem}}

A further consequence is

{\begin{theorem}  \cite{Ka2007} { Given a number $d$ and an algebraically closed field $K$, the class of $d$-dimensional tame strongly simply connected $K$-algebras induces a Zariski-open subset in the variety of associative $d$-dimensional $K$-algebras.}\end{theorem}}

Let us finish with presenting some general idea of connecting certain properties which seem to be of different origin at the first sight. In the above presentation we tried to highlight  connections between some general properties. Namely,  we know that:

\begin{enumerate}
\item Finite representation type is open \cite{Ga1974}.
\item Finite representation type is finitely axiomatizable \cite{JeLe1982}.
\item There is a criterion of representation-finiteness (expressed in terms of the Tits quadratic form) valid for  simply connected algebras.
\item Galois coverings behave nice with respect to representation-finiteness.
\end{enumerate}

Some time ago Jose Antonio de la Pe\~na asked  what remains of that in the class of  tame algebras? Let us write down that:

\begin{enumerate}
\item We can formulate the conjecture that tame is open.
\item This conjecture is equivalent to the conjecture that tame is finitely axiomatizable \cite{Ka2002}.
\item We have quadratic form criterion for tameness of strongly simply connected algebras - Theorem \ref{main}.
\item We can formulate the "Galois Covering Preserve Tameness (GCPT)" Conjecture: Assume that a torsion-free group $G$ acts freely on the objects of a locally bounded category $\widetilde{A}$. Then $A=\widetilde{A}/G$ is tame if and only if $\widetilde{A}$ is tame.
\end{enumerate}

Let us remark that
the "only if" part of GCPT is  known by \cite{DoSk1985}. The "if" has been  announced by Drozd and Ovsienko. But we are still - excluding Theorem \ref{main} - on a conjectural level. However there is an unexpected connection between some of the conjectures.

It is known that we can express the theory of Galois coverings in the language of graded algebras  \cite{Gr1983}.  This language is more suitable to our purposes.

Let us consider graded algebras over fields of characteristic $p$ (0 or a prime). We have the following

{\begin{theorem}\cite{KaPe2005} {  If tame graded algebras of fixed dimension form a finitely axiomatizable class, then GCPT holds for $p'$-residually finite group $G$.}\end{theorem}}

We say that a group $G$ is {\em $p'$-residually finite} if for any $a\in G$, $a\neq e$, there exists $H\vartriangleleft G$ of finite index not divisible by $p$ such that $a\notin H$.

{\sc Faculty of Mathematics and Computer Science, Nicolaus Copernicus University, Toru\'n, Poland}

skasjan@mat.umk.pl, skowron@mat.umk.pl


\begin{thebibliography}{99}
\bibitem{Au1978}  M. Auslander, Applications of morphisms determined by modules. Representation Theory of Algebras (Philadelphia, Pa., 1976),  245–327, Lecture Notes in Pure Appl. Math., Vol. 37, Dekker, New York, 1978.


\bibitem{AsSiSk2006} I. Assem, D. Simson and A. Skowroński, {\em Elements of the Representation Theory of Associative Algebras. Vol. 1: Techniques of Representation Theory.} London Mathematical Society Student Texts, Vol. 65. Cambridge University Press, Cambridge, 2006.

\bibitem{AsSk1988} I. Assem and A. Skowroński,  On some classes of simply connected algebras, Proc. London Math. Soc. 56 (3) (1988), no. 3, 417–450.

\bibitem{Ba1985}  R. Bautista,  On algebras of strongly unbounded representation type. Comment. Math. Helv. 60 (1985), no. 3, 392–399.

\bibitem{BaGaRoSa1985} R. Bautista, P. Gabriel, A. V. Roiter, and L. Salmerón,  Representation-finite algebras and multiplicative bases. Invent. Math. 81 (1985), no. 2, 217–285.

\bibitem{Bo1983}  K. Bongartz,  Algebras and quadratic forms. J. London Math. Soc. (2) 28 (1983), no. 3, 461–469.

\bibitem{Bo1984} K. Bongartz, A criterion for finite representation type. Math. Ann. 269 (1984), no. 1, 1–12.

\bibitem{Bo1985}  K. Bongartz,  Indecomposables are standard. Comment. Math. Helv. 60 (1985), no. 3, 400–410.

\bibitem{BoGa1981} K. Bongartz aand P. Gabriel,   Covering spaces in representation-theory. Invent. Math. 65 (1981/82), no. 3, 331–378.

\bibitem{Br1975}  S. Brenner, Quivers with commutativity conditions and some phenomenology of forms. Representations of Algebras (Ottawa, Ont., 1974), 29-53,   Lecture Notes in Math., Vol.  488, Springer, Berlin, 1975.


\bibitem{BrGa1984}  O. Bretscher and P. Gabriel,  The standard form of a representation-finite algebra. Bull. Soc. Math. France 111 (1983), no. 1, 21–40.

\bibitem{BrTo1986}  O. Bretscher and G. Todorov,  On a theorem of Nazarova and Roiter. Representation Theory I (Ottawa, Ont., 1984), 50–54, Lecture Notes in Math., Vol. 1177, Springer, Berlin, 1986.


\bibitem{BrPeSk2011} T. Brüstle, J. A. de la Pena, A. Skowroński,  Tame algebras and Tits quadratic forms. Adv. Math. 226 (2011), no. 1, 887–951.

\bibitem{BuRi1987}  M. C. R. Butler and C. M. Ringel,  Auslander-Reiten sequences with few middle terms and applications to string algebras. Comm. Algebra 15 (1987), no. 1-2, 145–179.

 \bibitem{CB1988} W. W. Crawley-Boevey,  On tame algebras and bocses. Proc. London Math. Soc. (3) 56 (1988), no. 3, 451–483.

\bibitem{CB1991} W. W.   Crawley-Boevey,  Tame algebras and generic modules. Proc. London Math. Soc. (3) 63 (1991), no. 2, 241–265.


 \bibitem{DoSk1985}  P. Dowbor and A. Skowroński,  On Galois coverings of tame algebras. Arch. Math. (Basel) 44 (1985), no. 6, 522–529.

 \bibitem{DoSk1987} P. Dowbor and A. Skowroński, Galois coverings of representation-infinite algebras. Comment. Math. Helv. 62 (1987), no. 2, 311–337.

 \bibitem{Dr1972}  Y. A. Drozd,  Matrix problems, and categories of matrices. (Russian) Investigations on the theory of representations. Zap. Naučn. Sem. Leningrad. Otdel. Mat. Inst. Steklov. (LOMI) 28 (1972), 144–153.

\bibitem{Dr1979}  Y. A. Drozd, Tame and wild matrix problems. (Russian) Representations and quadratic forms (Russian), pp. 39–74, 154, Akad. Nauk Ukrain. SSR, Inst. Mat., Kiev, 1979.

\bibitem{Fi1985}  U. Fischbacher, Une nouvelle preuve d'un théoreme de Nazarova et Roiter. (French) [A new proof of a theorem of Nazarova and Roiter] C. R. Acad. Sci. Paris Sér. I. Math. 300 (1985), no. 9, 259–262.

\bibitem{Ga1974}  P. Gabriel,  Finite representation type is open. Representations of Algebras (Ottawa, Ont., 1974), 132-155,    Lecture Notes in Math., Vol. 488, Springer, Berlin, 1975.



\bibitem{Ge1995} Ch. Geiss, On degenerations of tame and wild algebras. Arch. Math. (Basel) 64 (1995), no. 1, 11–16.

\bibitem{Gr1983}  E. L. Green,  Graphs with relations, coverings and group-graded algebras. Trans. Amer. Math. Soc. 279 (1983), no. 1, 297–310.

\bibitem{GrPr2016}  L. Gregory and M. Prest, Representation embeddings, interpretation functors and controlled wild algebras. J. Lond. Math. Soc. (2) 94 (2016), no. 3, 747–766.

\bibitem{HaPr2015}  R. Harland and M. Prest,  Modules with irrational slope over tubular algebras. Proc. Lond. Math. Soc. (3) 110 (2015), no. 3, 695–720.



\bibitem{JeLe1981}  Ch. U. Jensen and H. Lenzing, Applications of model theory to representations of finite-dimensional algebras. Math. Z. 178 (1981), no. 1, 83–98.

\bibitem{JeLe1982}  Ch. U. Jensen and H. Lenzing, Homological dimension and representation type of algebras under base field extension. Manuscripta Math. 39 (1982), no. 1, 1–13.

\bibitem{JeLe1989}  Ch. U. Jensen and H. Lenzing,  {\em Model-theoretic Algebra with Particular Emphasis on Fields, Rings, Modules}. Algebra, Logic and Applications, Vol. 2, Gordon and Breach Science Publishers, New York, 1989.

\bibitem{Ka2002}  S. Kasjan, On the problem of axiomatization of tame representation type. Fund. Math. 171 (2002), no. 1, 53–67.

\bibitem{Ka2007} S. Kasjan, Tame strongly simply connected algebras form an open scheme. J. Pure Appl. Algebra 208 (2007), no. 2, 435–443.


\bibitem{KaPa2014} S. Kasjan and G. Pastuszak, On the existence of super-decomposable pure-injective modules over strongly simply connected algebras of non-polynomial growth. Colloq. Math. 136 (2014), no. 2, 179–220.

\bibitem{KaPa2016} S. Kasjan and G. Pastuszak, Super-decomposable pure-injective modules over algebras with strongly simply connected Galois coverings. J. Pure Appl. Algebra 220 (2016), no. 8, 2985–2999.

\bibitem{KaPe2005} S. Kasjan and J. A. de la Pe\~na, Galois coverings and the problem of axiomatization of the representation type of algebras. Extracta Math. 20 (2005), no. 2, 137–150.

\bibitem{Ke1989} O. Kerner, Tilting wild algebras. J. London Math. Soc. (2) 39 (1989), no. 1, 29–47.

\bibitem{MaPe1983}  R. Martínez-Villa and J. A. de la Pe\~na,  The universal cover of a quiver with relations. J. Pure Appl. Algebra 30 (1983), no. 3, 277–292.

\bibitem{Na1973} L. A. Nazarova,  Representations of quivers of infinite type. (Russian) Izv. Akad. Nauk SSSR Ser. Mat. 37 (1973), 752–791.

\bibitem{NaRo1973} L. A.  Nazarova and A. V Roiter, Categorial matrix problems, and the Brauer-Thrall problem. (Russian)  Izdat. "Naukova Dumka'', Kiev, 1973. 100 pp., German translation  Mitt. Math. Sem. Giessen Heft 115 (1975).

\bibitem{NoSk1996}  R. N\"orenberg and A. Skowroński, Tame minimal non-polynomial growth strongly simply connected algebras. Colloq. Math. 73 (1997), no. 2, 301-330.

\bibitem{Pa2015}  G. Pastuszak,  Strongly simply connected algebras with super-decomposable pure-injective modules. J. Pure Appl. Algebra 219 (2015), no. 8, 3314–3321.

\bibitem{Pe1991} J. A. de la Pe\~na,  On the dimension of the module-varieties of tame and wild algebras. Comm. Algebra 19 (1991), no. 6, 1795–1807.


\bibitem{PeSk1996}  J. A. de la Pe\~na and A. Skowroński,  Geometric and homological characterizations of polynomial growth strongly simply connected algebras. Invent. Math. 126 (1996), no. 2, 287–296.

\bibitem{PeSk1999}  J. A. de la Pe\~na and A. Skowroński,  The Tits and Euler forms of a tame algebra. Math. Ann. 315 (1999), no. 1, 37–59.

\bibitem{PeSk2011}  J. A. de la Pe\~na and A. Skowroński,  The Tits forms of tame algebras and their roots. Representations of Algebras and Related Topics, 445–499, EMS Ser. Congr. Rep., Eur. Math. Soc., Zürich, 2011.

\bibitem{Pu2004}  G. Puninski,  Superdecomposable pure-injective modules exist over some string algebras. Proc. Amer. Math. Soc. 132 (2004), no. 7, 1891–1898.

\bibitem{Ri1979} C. M. Ringel,
On algorithms for solving vector space problems. I. Report on the Brauer-Thrall conjectures: Rojter's theorem and the theorem of Nazarova and Rojter. Representation Theory I (Ottawa, Ont., 1979),  104–136,
Lecture Notes in Math., Vol. 831, Springer, Berlin, 1980.

\bibitem{Ri1980} C. M. Ringel, Tame Algebras. Representation Theory, I (Ottawa, Ont., 1979),  137–287,
Lecture Notes in Math., 831, Springer, Berlin, 1980.



\bibitem{Ri1984} C. M. Ringel, Tame Algebras and Integral Quadratic Forms, Lecture Notes in Math. Vol. 1099, Springer, Berlin, Heidelberg, 1984.

\bibitem{Ro1968}  A. V. Roiter,  Unboundedness of the dimensions of the indecomposable representations of an algebra which has infinitely many indecomposable representations. (Russian) Izv. Akad. Nauk SSSR Ser. Mat. 32, 1968.

\bibitem{Ro1980}  A. V. Roiter,  Matrix problems and representations of BOCSs. Representation Theory I (Ottawa, Ont., 1979),  288-324,
Lecture Notes in Math., Vol. 831, Springer, Berlin, 1980.

\bibitem{Si1992} D. Simson, \emph{Linear Representations of Partially Ordered Sets and Vector Space Categories}.  Algebra, Logic and Applications, Vol. 4, Gordon \& Breach Science Publishers, 1992.

\bibitem{SiSk2007}   D. Simson and A. Skowro\'nski, {\em Elements of the Representation Theory of Associative Algebras, Vol. 3: Representation-Infinite Tilted Algebras}, London Mathematical Society Student Texts; Vol. 72, Cambridge University Press, 2007.

\bibitem{Sk1987} A. Skowroński, Group algebras of polynomial growth. Manuscripta Math. 59 (1987), no. 4, 499–516.

\bibitem{Sk1990} A. Skowroński, Algebras of polynomial growth. Topics in Algebra, Part 1 (Warsaw, 1988), 535–568, Banach Center Publ., Vol. 26, Part 1, PWN, Warsaw, 1990.

 \bibitem{Sk1992} A. Skowroński, Simply connected algebras and Hochschild cohomologies. in: Representations of Algebras (Ottawa, Ont., 1992), 431–447, CMS Conf. Proc., Vol. 14, Amer. Math. Soc., Providence, RI, 1993.


 \bibitem{Sk1997} A. Skowroński, Simply connected algebras of polynomial growth. Compositio Math. 109 (1997), no. 1, 99–133.

 \bibitem{Sk1998} A. Skowroński, Tame quasi-tilted algebras. J. Algebra 203 (1998), no. 2, 470–490.

\bibitem{Sm1980}  S. O. Smal\o,  The inductive step of the second Brauer-Thrall conjecture. Canad. J. Math. 32 (1980), no. 2, 342–349.

 \bibitem{WaWa1985}  B. Wald and J. Waschbüsch, Tame biserial algebras. J. Algebra 95 (1985), no. 2, 480–500.
\end{thebibliography}
\end{document}